\documentclass[leqno]{macrorend}

\volumeyear{xxxx}\yearnumber{x}\volumenumber{xx}

\usepackage[latin1]{inputenc}

\newcommand{\cI}{{\mathcal I}}

\newcommand{\cO}{{\mathcal O}}
\newcommand{\cE}{{\mathcal E}}
\newcommand{\cF}{{\mathcal F}}
\newcommand{\cM}{{\mathcal M}}
\newcommand{\cN}{{\mathcal N}}
\newcommand{\cL}{{\mathcal L}}

\begin{document}

%
\selectlanguage{english}

%
\articolo[Liaison with Cohen--Macaulay Modules]{Liaison with Cohen--Macaulay Modules}{Robin Hartshorne}

%

\begin{abstract}
We describe some recent work concerning Gorenstein liaison of codimension two subschemes of a projective variety.  Applications make use of the algebraic theory of maximal Cohen--Macaulay modules, which we review in an Appendix.
\end{abstract}

\section*{Introduction}

The purpose of this paper is to report on some recent work in the area of Gorenstein liaison.  For me this is a pleasant topic, because it illustrates the field of algebraic geometry at its best.  After all, algebraic geometry could be described as the use of algebraic techniques in geometry and the use of geometric methods to understand algebra.  In the work I describe here, we found an unexpected connection between the theory of maximal Cohen--Macaulay modules about which there is considerable algebraic literature, and the notion of Gorenstein liaison, which has emerged recently as geometers attempted to generalize results about curves in $\PP^3$ to varieties of higher codimension.

In \S \ref{sec1}, we review the ``classical'' case of curves in $\PP^3$.  In \S \ref{sec2} we describe generalizations of the notion of liaison to schemes of higher dimension and higher codimension.  Sections \ref{sec3} and \ref{sec4} develop the main new idea, which is instead of working directly with schemes of codimension $\ge 3$ in $\PP^n$, to consider subschemes of codimension $2$ of an arithmetically Gorenstein scheme $X$ in $\PP^n$.  Any liaison in $X$ is also a liaison in $\PP^n$, so this method is useful to establish existence of liaisons in $\PP^n$, but it cannot give negative results.  We hope that the study of liaison on $X$ may be interesting in its own right, and give more insight into the nature of liaison in general.

Section \ref{sec5} gives some applications, and \S \ref{sec6} describes an interesting open problem.  The algebraic theory of maximal Cohen--Macaulay modules is reviewed in an Appendix.

The principal new results described here are joint work with Marta Casanellas and Elena Drozd, given in detail in the papers \cite{CH} and \cite{CDH}.  For background on liaison, I recommend the book of Migliore \cite{M}, and for information on Cohen--Macaulay modules, the book of Yoshino \cite{Y}.

It was a pleasure to attend the conference Syzygy 2005 in Torino in honor of Paolo Valabrega's sixtieth birthday, and I dedicate this paper respectively to him.

\section{Curves in $\PP^3$}
\label{sec1}

We review the case of curves in $\PP^3$, which has been known for some time, as a model for the more general situations that we will consider below.

We work over an algebraically closed field $k$.  A {\em curve} is a purely one-dimensional scheme without embedded points.  If $C_1$ and $C_2$ are curves in $\PP^3$, we say they are {\em linked} by a complete intersection curve $Y$, if $C_1 \cup C_2 = Y$ and $\cI_{C_i,Y} \cong {\mathcal H}om(\cO_{C_j},\cO_Y)$ for $i,j = 1,2$, $i \ne j$.  The equivalence relation generated by chains of linkages is called {\em liaison}.  If a liaison is accomplished by an even number of linkages, it is called {\em even liaison}.

To any curve $C$ in $\PP^3$ we associate its {\em Rao module} $M_C = H_*^1(\cI_C) = \bigoplus_{n \in \ZZ} H^1(\PP^3,\cI_{C,\PP^3}(n))$.

The basic results about curves in $\PP^3$ are the following

\begin{theorem}[Rao \cite{R}]
\label{th1.1}
{\em a)} Two curves $C_1$ and $C_2$ are in the same even liaison equivalence class if and only if their Rao modules $M_1,M_2$ are isomorphic, up to a shift in degrees.

{\em b)} For any finite-length graded $R = k[x_0,x_1,x_2,x_3]$-module $M$, there exists a nonsingular irreducible curve $C$ in $\PP^3$, whose Rao module is isomorphic to a shift of $M$.

Thus the even liaison equivalence classes of curves in $\PP^3$ are in one-to-one correspondence with finite length graded $R$-modules, up to shift.
\end{theorem}

For the next statement we need the notion of biliaison.  If a curve $C_1$ lies on a surface $S$, and if $C_2 \sim C_1 + mH$, meaning linear equivalence in the sense of generalized divisors \cite{GD} on $S$, where $H$ is the hyperplane section of $S$, and $m$ is an  integer, then we say that $C_2$ is obtained by an {\em elementary biliaison} of height $m$ from $C_1$.  If $m \ge 0$ it is an {\em ascending} elementary biliaison.  It is easy to see that an elementary biliaison gives an even liaison between $C_1$ and $C_2$.  It is also easy to calculate numerical invariants of $C_2$, such as degree, genus, and postulation, from those of $C_1$ in terms of $m$ and the degree of $S$.

\begin{theorem}[Lazarsfeld--Rao property \cite{BBM}, \cite{MDP}]
\label{th1.2}
{\em a)} In any even liaison equivalence class of curves in $\PP^3$, the {\em minimal curves} (meaning those of minimal degree) form an irreducible family.

{\em b)} Any curve that is not minimal in its even liaison equivalence class can be obtained by a sequence of ascending elementary biliaisons from some minimal curve.
\end{theorem}

\begin{remark}
\label{rem1.3}
These results generalize well to subschemes $V$ of codimension two in $\PP^n$.  The Rao module has to be replaced by a series of higher deficiency modules $H_*^i(\cI_V)$ for $0 < i \le \dim V$ and certain extensions between them:  the best way to express this is by an element of the derived category.  Or one can use the so-called $\cE$-type resolution, in which case the set of even liaison equivalence classes of schemes $V$ of codimension two is in one-to-one correspondence with coherent sheaves $\cE$ (satisfying some additional conditions), up to stable equivalence and shift, and this in turn is in one-to-one correspondence with the quasi-isomorphism classes of certain complexes in the derived category replacing the Rao module.

The Lazarsfeld--Rao property also generalizes to codimension two subschemes of quite general schemes.  See for example \cite{HR} for precise statements and further references.
\end{remark}

\section{Generalizations}
\label{sec2}

When we consider curves in $\PP^4$, or more generally, subschemes of codimension $\ge 3$ in any $\PP^n$, the direct analogue of Rao's theorem fails.  There are infinitely many distinct even liaison equivalence classes of curves all having the same Rao module, which can be distinguished by other cohomological invariants \cite{KMMNP}.  It seems that liaison using complete intersections, as we have defined it, is much too rigid to give an analogous theory in higher codimension.

The notion of Gorenstein liaison seems to be a better candidate for generalizing the theory.

\begin{definition}
Two subschemes $V_1,V_2$ of $\PP^n$, equidimensional and without embedded components, are {\em $G$-linked} by an arithmetically Gorenstein scheme $Y$ (meaning the homogeneous coordinate ring of $Y$ is a Gorenstein ring) if $V_1 \cup V_2 = Y$ and $\cI_{V_i,Y} \cong {\mathcal H}om(\cO_{V_j},\cO_Y)$ for $i,j = 1,2$, $i \ne j$.  The equivalence relation generated by chains of $G$-links is called {\em Gorenstein liaison} (or $G$-liaison for short), and if a $G$-liaison can be accomplished by an even number of $G$-links, it is called {\em even $G$-liaison}.
\end{definition}

It is easy to see for curves in $\PP^n$ that even $G$-liaison preserves the Rao module (up to shift), as in the case of $\PP^3$, and this naturally leads to the converse problem:

\begin{problem}
\label{prob2.1}
If two curves in $\PP^n$ have isomorphic Rao modules (up to shift), are they in the same even $G$-liaison class?
\end{problem}

This problem is open at present.  The special case when the Rao module is zero is the case of {\em arithmetically Cohen--Macaulay} (ACM) curves, meaning that the homogeneous coordinate ring is a Cohen--Macaulay ring.  This includes in particular the complete intersection curves.  So the problem, which now can be stated for schemes of any dimension is

\begin{problem}
\label{prob2.2}
If $V$ is an ACM scheme in $\PP^n$, is $V$ in the Gorenstein liaison class of a complete intersection (glicci for short)?
\end{problem}

This problem is also open at present, though many special cases are known (see for example \cite{KMMNP}).  There are also candidates for counterexamples (as yet unproven), such as 20 general points in $\PP^3$, or a general curve of degree 20 and genus 26 in $\PP^4$ \cite{SEG}.

Our approach in this paper, instead of studying the problem directly in $\PP^n$, will be to study codimension two subscheme of an arithmetically Gorenstein variety $X$ in $\PP^n$.  Liaisons in $X$ can also be considered to be liaisons in $\PP^n$, and thus we study the problem of higher codimension subschemes in $\PP^n$ indirectly.  While most of our results are valid for $X$ of any dimension, for simplicity in this paper we will stick to dimension $3$.

So here is the set-up.  Let $X$ be a fixed normal arithmetically Gorenstein subvariety of dimension $3$ in $\PP^n$.  We also keep fixed the embedding and hence the sheaf $\cO_X(1)$ on $X$ that defines the class of a hyperplane section $H$ of $X$.

If $C_1$ and $C_2$ are curves in $X$, we say that $C_1$ and $C_2$ are {\em linked} by a curve $Y$ in $X$ if $C_1 \cup C_2 = Y$ and $\cI_{C_i,Y} \cong {\mathcal H}om(\cO_{C_j},\cO_Y)$ for $i,j = 1,2$, $i \ne j$.  If $Y$ is a {\em complete intersection} in $X$, meaning that $Y$ is the intersection of surfaces defined by sections of $\cO_X(a)$, $\cO_X(b)$ in $X$, then we say it is a {\em $CI$-linkage}.  If $Y$ is arithmetically Gorenstein (in the ambient $\PP^n$), it is a {\em $G$-linkage}.  These linkages give rise to the equivalence relations of {\em $CI$-liaison} and {\em even $CI$-liaison} and {\em $G$-liaison} and {\em even $G$-liaison} as before.

Note that a $CI$-liaison in $X$ is not necessarily a $CI$-liaison in $\PP^n$, unless $X$ itself is a complete intersection.  However, a $G$-liaison in $X$ is also a $G$-liaison in $\PP^n$.

If $S$ is a surface in $X$ containing a curve $C$, and if $C'$ is another curve on $S$, with $C' \sim C + mH$, meaning linear equivalence of generalized divisors on $S$, where $H$ is the hyperplane section, we say $C'$ is obtained from $C$ by an {\em elementary biliaison} from $C$.  If $S$ is a complete intersection in $X$ (corresponding to $\cO_X(a)$ for some $a$) it is a {\em $CI$-biliaison}.  If $S$ is an ACM scheme (in $\PP^n$) it is a {\em $G$-biliaison}.  It is easy to see that a $CI$-biliaison is an even $CI$-liaison.  In fact, the equivalence relation generated by $CI$-biliaisons is the same as even $CI$-liaison (proof similar to \cite[$4.4$]{GD}).  One can show also that a $G$-biliaison is an even $G$-liaison \cite{KMMNP}, \cite[$3.6$]{GDB}, however in general the equivalence relation generated by $G$-biliaisons is not the same as even $G$-liaison, as we can see from the following example.

\begin{example}
\label{exam2.3}
Let $X$ be a nonsingular quadric hypersurface in $\PP^4$.  Every surface on $X$ is a complete intersection, and in particular has even degree.  Thus $G$-biliaisons preserve the parity of the degree of a curve.  On the other hand, the union of a rational quartic curve with a line meeting it at two points is an arithmetically Gorenstein elliptic quintic, so the two curves are $G$-linked.  One line can also be linked to another line by a conic, so we see that even $G$-liaison does not preserve parity of degree.
\end{example}

In studying $G$-liaison and $G$-biliaison on $X$, an important role is played by the category of ACM sheaves on $X$.  An ACM {\em sheaf} is a coherent sheaf $\cE$ on $X$ that is locally Cohen--Macaulay and has vanishing intermediate cohomology:  $H_*^i(\cE) = 0$ for $i = 1,2$.  If $X$ is $\PP^3$, the only ACM sheaves are the {\em dissoci\'e} sheaves, i.e., direct sums of line bundles $\cO_X(a_i)$, by a theorem of Horrocks.  However, if $X$ is not $\PP^3$, there are others, and the category of these sheaves reflects interesting properties of $X$.

To see why these sheaves are important for $G$-liaison and $G$-biliaison, we first mention the following result relating them to ACM surfaces in $X$ and arithmetically Gorenstein (AG) curves in $X$.

\begin{proposition}
\label{prop2.4}
{\em a)} If $S$ is an {\em ACM} surface in $X$, then its ideal sheaf $\cI_{S,X}$ is a rank $1$ {\em ACM} sheaf on $X$.  Conversely if $\cL$ is a rank $1$ {\em ACM} sheaf on $X$, then for any $a \gg 0$, the sheaf $\cL(-a)$ is isomorphic to the ideal sheaf $\cI_{S,X}$ of an {\em ACM} surface in $X$.

{\em b)} If $Y$ is an {\em AG} curve in $X$, then there is an exact sequence
\[
0 \to \cO_X(-a) \to \cN \to \cI_{Y,X} \to 0
\]
for some $a \in \ZZ$, where $\cN$ is a rank $2$ {\em ACM} sheaf on $X$ with $c_1(\cN) = -a$.  Conversely if $\cN$ is any {\em orientable} (meaning $c_1(\cN) = \cO_X(a)$ for some $a \in \ZZ$) rank $2$ {\em ACM} sheaf on $X$, and if $s$ is a sufficiently general section of $\cN(a)$ for $a \gg 0$, then $s$ induces an exact sequence
\[
0 \to \cO_X \stackrel{s}{\to} \cN(a) \to \cI_{Y,X}(b) \to 0
\]
for some {\em AG} curve $Y$ in $X$ and some $b \in \ZZ$.
\end{proposition}

\begin{proof}
Part a) is elementary, while part b) is the usual Serre correspondence \cite[$2.9$]{CDH}.
\end{proof}

Thus we see that the ACM surfaces and AG curves, which are used to define $G$-biliaison and $G$-liaison, respectively, correspond in a natural way to rank $1$ and rank $2$ ACM sheaves on $X$.  In the following two sections, we will study $G$-biliaison and $G$-liaison separately.

\section{Gorenstein biliaison}
\label{sec3}

As in the previous section, we consider a normal arithmetically Gorenstein $3$-fold $X$, and we will consider Gorenstein biliaison of curves on $X$.

First of all, let's see what happens with a single elementary Gorenstein biliaison.  Let $S$ be an ACM surface in $X$, let $C$ be a curve in $S$, and let $C' \sim C + mH$ on $S$.  Then by construction, $\cI_{C',S} \cong \cI_{C,S}(-m)$.  Thus we can write exact sequences
\[
\begin{array}{ccccccccc}
0 &\to &\cI_S &\to &\cI_{C'} &\to &\cI_{C',S} &\to &0 \\
& & & & & &\| \\
0 &\to &\cI_S(-m) &\to &\cI_C(-m) &\to &\cI_{C,S}(-m) &\to &0.
\end{array}
\]
If we let $\cF$ be the fibered sum of $\cI_{C'}$ and $\cI_C(-m)$ over $\cI_{C',S} = \cI_{C,S}(-m)$, we obtain sequences
\[
\begin{array}{ccccccccc}
0 &\to &\cI_S &\to &\cF &\to &\cI_C(-m) &\to &0 \\
0 &\to &\cI_S(-m) &\to &\cF &\to &\cI_{C'} &\to &0.
\end{array}
\]
Note here that the same coherent sheaf $\cF$ appears in the middle of each sequence, and that the sheaves on the left are rank $1$ ACM sheaves on $X$ that are isomorphic, up to twist.

Conversely, given exact sequences
\[
\begin{array}{ccccccccc}
0 &\to &\cL &\to &\cF &\to &\cI_C(a) &\to &0 \\
0 &\to &\cL' &\to &\cF &\to &\cI_{C'}(a') &\to &0
\end{array}
\]
with the same coherent sheaf $\cF$ in the middle, where $C,C'$ are curves in $X,a,a'$ integers, and $\cL,\cL'$ rank $1$ ACM sheaves that are isomorphic up to twist, it follows that $C'$ is obtained by a single elementary $G$-biliaison from $C$.  The idea of proof is to consider the composed map $\cL' \to \cF \to\cI_C(a)$.  If this map is $0$, then $C' = C$, which is a trivial $G$-biliaison.  If it is not zero, composing with the inclusion $\cI_C(a) \subseteq \cO_X(a)$ identifies $\cL'(-a)$ with the ideal sheaf $\cI_S$ of an ACM surface on $X$ and then one sees easily that $C' \sim C + (a'-a)H$ on $S$ \cite[$3.1,3.3$]{CH}.

With a little more work, one can arrive at an analogous criterion for two curves to be related by a finite succession of elementary $G$-biliaisons.

\begin{theorem}[$\mbox{\cite[$3.1$]{CH}}$]
\label{th3.1}
Two curves $C,C'$ on the normal arithmetically Gorenstein $3$-fold $X$ are in the same Gorenstein biliaison equivalence class if and only if there exist exact sequences
\[
\begin{array}{ccccccccc}
0 &\to &\cE &\to &\cN &\to &\cI_C(a) &\to &0 \\
0 &\to &\cE' &\to &\cN &\to &\cI_{C'}(a') &\to &0
\end{array}
\]
with the same coherent sheaf $\cN$ in the middle, where $a,a'$ are 
integers, and where $\cE$ and $\cE'$ are {\em ACM} sheaves each having a filtration whose quotients are rank $1$ {\em ACM} sheaves (we call them {\em layered} {\em ACM} sheaves), and such that the rank $1$ quotients of these filtrations of $\cE$ and $\cE'$ are isomorphic up to order and twists.
\end{theorem}

The point is that each $G$-biliaison contributes a rank $1$ ACM factor, but that these make up the two sheaves $\cE$ and $\cE'$ in a different order, and with different twists associated to each.

If $\cE$ is a layered ACM sheaf as above, the filtration with rank $1$ ACM quotients may not be unique.  Taking advantage of this are two ``exchange lemmas'' \cite[$3.4,4.6$]{CH} that allow one to replace one $\cE$ by another $\cE'$ having the same factors, in sequences as in Theorem~\ref{th3.1}, after passing to another curve in the same $G$-biliaison class.  These form a sort of converse to Theorem~\ref{th3.1}, and allow us to formulate a necessary and sufficient condition for the property analogous to Problem~\ref{prob2.2} on $X$, namely that every ACM curve on $X$ should be in the $G$-biliaison class of a complete intersection on $X$.  This condition is a bit complicated to state (see \cite[$4.2,4.3$]{CH}), so instead here we will explain the result only in one interesting special case.

\begin{theorem}[$\mbox{\cite[$6.2$]{CH}}$]
\label{th3.2}
Let $X$ be the cone over a nonsingular quadric surface in $\PP^3$.  (Thus $X$ is a normal quadric hypersurface  in $\PP^4$ having one double point.)  Then two curves $C$ and $C'$ on $X$ are in the same Gorenstein biliaison equivalence class if and only if their Rao modules are isomorphic, up to shift.  In particular, all {\em ACM} curves are equivalent for $G$-biliaison.
\end{theorem}

\medskip
\noindent
{\em Idea of Proof.}  It is obvious that Gorenstein biliaison preserves the Rao module, up to shift, so one direction is clear.

For the other direction, let $C$ be any curve in $X$, with Rao module $M$.  Our strategy is to construct another curve $C'$ that depends only on $M$, and then show that $C$ and $C'$ are in the same $G$-biliaison equivalence class, which will prove the theorem.

Given $M$, let $M^*$ be the dual module, and take a resolution
\[
0 \to G \to F_2 \to F_1 \to F_0 \to M^* \to 0
\]
over $R$, the homogeneous coordinate ring of $X$, where the $F_i$ are free graded $R$-modules, and $G$ is the kernel.  Let $\cN'$ be the sheaf associated to $G^{\vee}$, and let $\cL'$ be a dissoci\'e sheaf of rank one less mapping to $\cN'$ so as to define a curve $C'$ by its cokernel:
\[
0 \to \cL' \to \cN' \to \cI_{C'}(a') \to 0.
\]

On the other hand, let
\[
0 \to \cL \to \cN \to \cI_C \to 0
\]
be an $\cN$-{\em type resolution} of $C$, i.e., with $\cL$ dissoci\'e and $\cN$ coherent, locally Cohen--Macaulay, and $H_*^1(\cN) \cong M$ and $H_*^2(\cN) = 0$.  Let $N = H_*^0(\cN)$ and take a resolution
\[
0 \to P \to L_1 \to L_0 \to N \to 0
\]
over $R$ with $L_i$ free and $P$ the kernel.  Dualizing gives an exact sequence
\[
0 \to N^{\vee} \to L_0^{\vee} \to L_1^{\vee} \to P^{\vee} \to M^* \to 0.
\]

Now there is a natural map of the earlier free resolution of $M^*$ into this one, and this gives us a map of $G$ to $N^{\vee}$, from which we obtain a natural map $\cN \to \cN'$.  By adding extra free factors if necessary, we may assume it is surjective, and then let $\cE$ be the kernel:
\[
0 \to \cE \to \cN \to \cN' \to 0.
\]
Since $\cN$ and $\cN'$ both have $H_*^1 = M$ and $H_*^2 = 0$, we see that $\cE$ is an ACM sheaf on $X$.  Furthermore, taking the composed map from $\cN$ to $\cI_{C'}(a')$ we obtain an exact sequence
\[
0 \to \cE \oplus \cL \to \cN \to \cI_{C'}(a') \to 0.
\]

In order to apply the criterion of Theorem~\ref{th3.1} we now need to use the special property of the quadric $3$-fold $X$ (see Appendix), which tells us first that every ACM sheaf on $X$ is layered, secondly that the only rank $1$ ACM sheaves on $X$ (up to twist) are $\cO_X$, $\cI_D$, and $\cI_E$, where $D,E$ represent the two types of planes in $X$, and thirdly that there is an exact sequence
\[
0 \to \cI_D \to \cO_X^2 \to \cI_E(1) \to 0.
\]

In the ACM sheaf $\cE$, copies of $\cI_D$ and $\cI_E$ (and their twists) must occur in equal numbers, because $\cE$ is orientable.  Then the exchange lemmas referred to above allow us to replace an $\cI_D$ plus an $\cI_E$ by an $\cO_X^2$.  Thus $\cE \oplus \cL$ is replaced by a dissoci\'e sheaf, and then Theorem~\ref{th3.1} tells us that $C$ and $C'$ are in the same $G$-biliaison class.  (For more details see \cite[$4.7,6.2$]{CH}.)

\section{Gorenstein liaison}
\label{sec4}

Let us consider a normal AG $3$-fold $X$, as before, and study Gorenstein liaison equivalence of curves in $X$.  Since the AG curves in $X$ are associated to rank $2$ ACM sheaves on $X$, as we saw above, we expect to see them play a role.

First of all, let us see what happens with a single Gorenstein liaison.  We track this behavior using the $\cN$-type resolution of a curve $C$.

\begin{proposition}
\label{prop4.1}
Let $C$ be a curve in $X$ with $\cN$-type resolution $0 \to \cL \to \cN \to \cI_C \to 0$, and suppose that $C$ is linked to a curve $C'$ by the {\em AG} curve $Y$.  Then $C'$ has an $\cN$-type resolution of the form
\[
0 \to \cL' \to \cN' \to \cI_{C'}(a') \to 0
\]
with $\cL'$ dissoci\'e, and where $\cN'$ is an extension
\[
0 \to \cL^{\vee} \oplus \cE^{\vee} \to \cN' \to \cN^{\sigma\vee} \to 0,
\]
where $\cE$ is the rank $2$ {\em ACM} sheaf associated to $Y$, and $\cN^{\sigma\vee}$ denotes the dual of the first syzygy sheaf of $\cN$.
\end{proposition}

To prove this (see \cite[$3.2$]{CDH}) one first uses the usual cone construction of the map $\cI_Y \subseteq \cI_C$, and this gives the sequence
\[
0 \to \cN^{\vee} \to \cL^{\vee} \oplus \cE^{\vee} \to \cI_{C'}(a) \to 0.
\]
This is not an $\cN$-type resolution, but by using the syzygy sheaf $\cN^{\sigma}$ of $\cN$, one can transform it into the desired $\cN$-type resolution.

Note what happens to the Rao module $M$.  From the definition of the syzygy sheaf
\[
0 \to \cN^{\sigma} \to \cF \to \cN \to 0
\]
with $\cF$ dissoci\'e, we see that $M \cong H_*^1(\cN) \cong H_*^2(\cN^{\sigma})$.  By Serre duality then $H_*^1(\cN^{\sigma\vee}) \cong M^*$, the dual of $M$, and this shows that the Rao module of $C'$  is $M^*$ shifted, as we would expect from a single liaison.

This proposition shows us that a single $G$-liaison complexifies the $\cN$-type resolution by throwing in a dual of a syzygy, and adding an extension by a rank $2$ ACM sheaf.  There is a sort of converse to this, showing how to simplify an $\cN$-type resolution by removing a rank $2$ ACM sheaf.  In general, this cannot be accomplished by a single $G$-liaison, but requires a more complicated procedure.

\begin{proposition}
\label{prop4.2}
Let $C$ be a curve with an $\cN$-type resolution
\[
0 \to \cL \to \cN \to \cI_C \to 0,
\]
and suppose given an exact sequence
\[
0 \to \cE \to \cN \to \cN' \to 0
\]
with $\cE$ a rank $2$ {\em ACM} sheaf and $\cN'$ a locally {\em CM} sheaf of rank $\ge 2$.  Then there is a curve $C'$ in the same even $G$-liaison equivalence class as $C$ having an $\cN$-type resolution
\[
0 \to \cL' \to \cN' \to \cI_{C'}(a') \to 0.
\]
\end{proposition}

\begin{proof}
See \cite[$3.4$]{CDH}.
\end{proof}

Using these two propositions, it is possible to give a criterion, in terms of the $\cN$-type resolutions, for when two curves are in the same $G$-liaison class \cite[$5.1$]{CDH}.  The exact statement, which involves successive extensions by rank $2$ ACM sheaves and their syzygy duals (which may no longer be of rank $2$), is rather complicated, so we omit it here.  Using this theorem, one can also give a criterion for every ACM curve to be in the Gorenstein liaison class of a complete intersection \cite[$5.4$]{CDH}.  Here we will just give one special case, albeit an interesting one.

\begin{theorem}
\label{th4.3}
Let $X$ be a nonsingular quadric hypersurface in $\PP^4$.  Then two curves are in the same even $G$-liaison class if and only if their Rao modules are isomorphic, up to shift.
\end{theorem}

\medskip
\noindent
{\em Sketch of Proof} (cf. \cite[$6.2$]{CDH}).  Let $C$ be any curve, with Rao module $M$, and let $C'$ be another curve with the same Rao module $M$, constructed as in the proof of Theorem~\ref{th3.2} above.  Following the plan of that proof we have $\cN$-type resolutions
\[
0 \to \cL \to \cN \to \cI_C \to 0
\]
\[
0 \to \cL' \to \cN' \to \cI_{C'}(a') \to 0
\]
and an exact sequence
\[
0 \to \cE \to \cN \to \cN' \to 0
\]
where $\cE$ is an ACM sheaf on $X$.

Now we invoke the special property of the nonsingular quadric $3$-fold $X$, which is that every ACM sheaf is a direct sum of a dissoci\'e sheaf and copies of twists of a single rank $2$ ACM sheaf $\cE_0$, associated to a line in $X$ (see Appendix).  We apply Proposition~\ref{prop4.2} repeatedly to remove copies of $\cE_0$ and its twists from $\cN$, thus eventually obtaining a curve $C''$, in the same even $G$-liaison class as $C$, and having an $\cN$-type resolution whose middle sheaf $\cN''$ differs from $\cN'$ only by a dissoci\'e sheaf.  Then $\cN'$ and $\cN''$ are stably equivalent, and so $C''$ and $C'$ are in the same even $CI$-liaison class, by Rao's theorem, and a fortiori in the same even $G$-liaison class.

Note that in the case of the nonsingular quadric $3$-fold, $G$-biliaison is just the same as $CI$-biliaison, hence is much too restrictive to provide a result like this theorem.

\section{Applications}
\label{sec5}

In \cite[$8.10$]{KMMNP} the authors showed, by an exhaustive listing of all possible ACM curves on these surfaces, that any ACM curve lying on a general smooth rational ACM surface in $\PP^4$ is glicci.  The rational ACM surfaces in $\PP^4$ (not counting those in $\PP^3$, for which the theorem is known) are the cubic scroll, the Del Pezzo surface of degree $4$, the Castelnuovo surface of degree $5$, and the Bordiga surface of degree $6$.

For the cubic scroll, the Del Pezzo, and the Castelnuovo surface, this result is an immediate consequence of our Theorem~\ref{th3.2}, because each of these surfaces is contained in a quadric $3$-fold with one double point.  Our method does not apply to the Bordiga surface, which is not contained in any quadric hypersurface.

\bigskip
In his paper \cite{L}, Lesperance studied curves in $\PP^4$ of the following form.  Let $C$ be the disjoint union $C_1 \cup C_2$ of two plane curves $C_1,C_2$, lying in two planes that meet at a single point $P$.  Let them have degrees $d_1,d_2$, and assume either a) $2 \le d_1 \le d_2$ or b) $2 \le d_1$ and $C_2$ contains the point $P$.  The Rao module is then $M \cong R/(I_P + R_{\ge d_1})$, which depends only on the point $P$ and the integer $d_1$.  Lesperance shows that all the curves of type a) and some of those of type b) are in the same Gorenstein liaison class, by using explicitly constructed $G$-liaisons.

Since a union of two planes meeting at a point is contained in a quadric hypersurface with one double point, it follows from our Theorem~\ref{th3.2} that all the above curves with the same Rao module are equivalent for $G$-liaison \cite[$6.4$]{CH}.

\bigskip
A third application is the following

\begin{theorem}
\label{th5.1}
Any arithmetically Gorenstein scheme $V$ in $\PP^n$ is in the Gorenstein liaison class of a complete intersection (glicci).
\end{theorem}

For the proof \cite[$7.1$]{CDH} we use the higher-dimensional analogues of the results described in this paper for an AG $3$-fold $X$.  By a Bertini-type theorem of Altman and Kleiman, one can find a complete intersection scheme $X$ in $\PP^n$, containing $V$, of dimension two greater than $V$, and smooth outside of $V$.  Then $X$ is normal and AG, and $V$ is a codimension two AG scheme in $X$, so there is an exact sequence
\[
0 \to \cO_X(-a) \to \cE \to \cI_{V,X} \to 0
\]
where $\cE$ is a rank $2$ ACM sheaf on $X$.  Let $\cM$ be a rank $2$ dissoci\'e sheaf on $X$ and consider the new $\cN$-type resolution of $V$,
\[
0 \to \cO_X(-a) \oplus \cM \to \cE \oplus \cM \to \cI_{V,X} \to 0.
\]

Then we apply the analogue of Proposition~\ref{prop4.2}, which is \cite[$3.4$]{CDH}, to remove $\cE$ and obtain another subscheme $V' \subseteq X$, in the same even $G$-liaison class as $V$, with an $\cN$-type resolution
\[
0 \to \cL' \to \cM \to \cI_{V',X}(a') \to 0.
\]
Since $\cM$ is rank $2$ dissoci\'e, it follows that $V'$ is a complete intersection in $X$, and since $X$ is itself a complete intersection in $\PP^n$, $V'$ is also a complete intersection in $\PP^n$.  Since $G$-liaisons in $X$ are also $G$-liaisons in $\PP^n$, we find that $V$ is glicci, as required.

\section{An open problem}
\label{sec6}

If there is a moral to all  the investigations of Gorenstein liaison so far, it seems to me that good results are obtained for schemes with some special structure, such as determinantal schemes \cite{KMMNP}, or schemes of codimension $2$ in low-degree hypersurfaces, such as the ones considered in \S \ref{sec3},\ref{sec4} above.

To describe a situation on the border between what is known and what is not known, I would like to consider the case of zero-dimensional subschemes of a non-singular cubic surface in $\PP^3$.  Though of one dimension lower than the discussions earlier in this paper, I think it is a good arena to test the essential difficulties of the subject.

So, let $X$ be a nonsingular cubic surface in $\PP^3$.  We consider zero-schemes $Z \subseteq X$.  Any zero-scheme is ACM, so there are two problems to consider.

\begin{problem}
\label{prob6.1}
Is every zero-scheme $Z \subseteq X$ in the $G$-biliaison equivalence class of a point?
\end{problem}

\begin{problem}
\label{prob6.2}
Is every zero-scheme $Z \subseteq X$ in the $G$-liaison equivalence class of a point?
\end{problem}

Both problems are open at present.  I will discuss what is known about them so far.

Using explicit $G$-liaisons and $G$-biliaisons on ACM curves on $X$, one can show that any set $Z$ of $n$ points in general position on $X$ is $G$-liaison equivalent to a point \cite[$2.4$]{SEG}.  The proof of this result is curious, in that one uses sequences of liaisons where the number of points may have to increase before it decreases.  For example, starting with 18 general points, one makes links to the following numbers of points (always in general position):  $18 \to 20 \to 28 \to 22 \to 16 \to 13 \to 7 \to 5 \to 3 \to 1$.  For points in special position, it seems hopeless to generalize this method.

Another approach, more in the spirit of this paper, is to study the category of ACM sheaves on $X$.  Faenzi \cite{F} has classified the rank $2$ ACM sheaves on $X$.  Up to twist, there is a finite number of possible Chern classes, and for fixed Chern classes, the possible sheaves form algebraic families of dimensions $\le 5$.  Already the presence of families of dimension $> 1$ shows that we are in a situation of ``wild CM-type'' (see Appendix).  Looking at Faenzi's results, again it seems hopeless to achieve a complete classification of ACM sheaves of all ranks on $X$.  One can show, however, that there are families of arbitrarily high dimension of indecomposable ACM sheaves of higher rank.

However, to answer the two problems above, one would not need a complete classification of ACM sheaves on $X$.  For an affirmative answer to Problem~\ref{prob6.1}, it would be sufficient to show \cite[$4.3$]{CH}.

$(*)$ Every orientable ACM sheaf $\cE$ on $X$ has a resolution
\[
0 \to \cF_2 \to \cF_1 \to \cE \to 0
\]
where $\cF_1$ and $\cF_2$ are layered ACM sheaves (i.e., successive extensions of rank $1$ ACM sheaves).

In regard to this property, there are examples of rank $2$ ACM sheaves $\cE$ on $X$, that are not layered themselves, but do have a resolution of this form.  So there seems to be some hope that this may hold.

For an affirmative answer to Problem~\ref{prob6.2}, it would be sufficient to show \cite[$5.4$]{CDH}.

$(**)$ Every orientable ACM sheaf on $\cE$ is stably equivalent to a {\em double-layered} sheaf on $X$ (which is a successive extension of rank $2$ ACM sheaves and their syzygies).

\section*{Appendix.  MCM modules and ACM sheaves}
\label{appen}

In this appendix we give a brief outline of some algebraic results that are needed to justify the results on ACM sheaves on quadric hypersurfaces used in \S \ref{sec3},\ref{sec4} above.

Let $R,{\mathfrak m}$ be a Cohen--Macaulay local ring.  A {\em maximal Cohen--Macaulay} (MCM) {\em module} is a finitely generated $R$-module $M$, with depth $M = \dim R$ and $\mbox{Supp } M = \mbox{Spec } R$.

For example, if $R$ is a regular local ring, every MCM module has homological dimension zero, and so is free.  Conversely, if $R,{\mathfrak m}$ is a Cohen--Macaulay local ring over which every MCM is free, then $R$ is regular.  Indeed, for any $R$-module $N$, consider a free resolution of length $n = \dim R$, and let $M$ be the kernel at the last step.  Then $M$ is an MCM module, hence free, and so $N$ has a finite free resolution.  Thus the ring $R$ has finite global homological dimension, and by a theorem of Serre, this implies that $R$ is regular.

Thus the presence of non-trivial MCM modules characterizes non-regular local rings, and the category of MCM modules is an interesting measure of the complexity of the singularity of the local ring.

In certain circumstances, Y.~Drozd \cite{D} has shown that local rings can be divided into three classes, depending on the behavior of the MCM modules.  It is true in any case that an MCM module can be written uniquely as a direct sum of {\em indecomposable} MCM modules, namely those that allow no further direct sum decomposition.  We say that $R$ is of {\em finite} CM-{\em type} if there is only a finite number of indecomposable MCM modules.  We say $R$ is of {\em tame} CM-{\em type} if the indecomposable MCM modules form a countable number of families of dimension at most one.  We say $R$ is of {\em wild} CM-{\em type} if there are families of arbitrarily large dimension of indecomposable MCM modules.  The tame-wild dichotomy theorem says (in certain cases) that only these cases can occur.  While to my knowledge this has not been proved in general, we can keep it in mind as a principle of what to expect when studying MCM modules.

The same definitions apply to the case of graded rings and graded modules, and thus admit a translation into sheaves on projective schemes.  If $X$ is a ACM scheme in $\PP^n$, we have defined an ACM {\em sheaf} on $X$ to be a locally Cohen--Macaulay coherent sheaf $\cE$ on $X$ with no intermediate cohomology:  $H_*^i(\cE) = 0$ for $0 < i < \dim X$.  If $R$ is the homogeneous coordinate ring of $X$, then we obtain a correspondence between ACM sheaves on $X$ and graded MCM modules on $R$ by sending a sheaf $\cE$ to the module $E = H_*^0(\cE)$, and sending the module $E$ to the associated sheaf $\cE = {\tilde E}$.  In carrying over results and definitions from the local case, we should consider graded modules up to shift, and ACM sheaves up to twist.  So we can say $X$ is of finite CM-type if there is only a finite number (up to twist) of isomorphism classes of indecomposable ACM sheaves.

To illustrate the different CM-types in the projective case, note that if $X$ is a nonsingular curve in $\PP^n$, then an ACM sheaf on $X$ is just a locally free sheaf $\cE$, also called a vector bundle.  If $X$ is rational, of degree $d$, the only indecomposable vector bundles (up to twist by $\cO_X(1)$) are line bundles of degrees $0 \le e < d$, so $X$ is of finite CM-type.  If $X$ is an elliptic curve, then by the classification theorem of Atiyah, for each rank and degree $(\mbox{mod } d = \mbox{deg } X)$ there is a one-parameter family of isomorphism classes of indecomposable vector bundles of rank $r$ and degree $e$.  Thus $X$ is of tame CM-type.  And if the genus of $X$ is $g \ge 2$, then as the rank grows, so does the dimension of the moduli space of stable vector bundles, so $X$ is of wild CM-type.

In the complex-analytic and complete local ring case, those local rings of isolated hypersurface singularities of finite CM type have been classified \cite{CMMI}, \cite{CMMII}.  They are the local rings of simple singularities in the sense of Arnol'd; in each dimension they are associated with Dynkin diagrams $A_n,D_n,E_6,E_7,E_8$, and their equations can be written explicitly.

Carrying these results over to the graded case, one obtains a list of all projective schemes of finite CM-type \cite{Y}, namely, projective spaces, nonsingular quadric hypersurfaces in any dimension, the rational cubic scroll in $\PP^4$, and the Veronese surface in $\PP^5$.  Furthermore, the indecomposable ACM sheaves on these varieties can be described explicitly, and this is where we find that there is just one non-trivial indecomposable ACM sheaf on the nonsingular quadric $3$-fold, mentioned in the proof of Theorem~\ref{th4.3}.

The main tools for studying MCM modules on hypersurface singularities, or ACM sheaves on hypersurfaces in $\PP^n$, are the matrix factorization, and the double branched covers and periodicity theorems of Kn\"orrer \cite{CMMI}.  We explain these in the projective case.

Let $X$ be a hypersurface in $\PP^n$, and let $\cE$ be an ACM sheaf on $X$.  Since the associated graded module $E = H_*^0(\cE)$ has depth $n$ over the coordinate ring $P = k[x_0,\dots,x_n]$ of $\PP^n$, there is a resolution
\[
0 \to \cL_1 \stackrel{\varphi}{\to} \cL_0 \to \cE \to 0
\]
by dissoci\'e sheaves $\cL_i$ on $\PP^n$ of the same rank $m$.  This gives a square matrix $\varphi$ of homogeneous forms in $P$.  Then one shows that there is another matrix $\psi$ of the same rank, with the property that $\psi \cdot \varphi = \varphi \cdot \psi = f \cdot id$, where $f$ is the equation of the hypersurface.  This is called a {\em matrix factorization} of $f$.  One sees also that $\mbox{det } \varphi = f^r$, where $r =$ rank $\cE$.  These constraints allow one to gain information about the possible ACM sheaves $\cE$ when the numbers are small enough.

The other technique is Kn\"orrer's double branched cover, and periodicity theorems, which allow one to pass from a hypersurface $X$ in $\PP^n$ defined by a polynomial $f \in P$ to the hypersurface $X'$ in $\PP^{n+1}$ defined by $f + x^2$, or the hypersurface $X''$ in $\PP^{n+2}$ defined by $f+x^2+y^2$, where $x$ and $y$ are new variables.

In the paper \cite{CH}, we use these techniques to show that the singular quadric $3$-fold $X$ in $\PP^4$ with one double point is of {\em countable} CM-{\em type}, namely it has only countably many indecomposable ACM sheaves (up to twist), and these are $\cO_X,\cI_D,\cI_E$, where $D,E$ are the two types of planes in $X$, and two infinite sequences $\cE_{\ell}$ and $\cE'_{\ell}$, for $\ell = 1,2,\dots$, of rank $2$ ACM sheaves that are each extensions of suitable twists of $\cI_D$ and $\cI_E$ \cite[$6.2$]{CH}, hence layered.  This is the result needed for the proof of Theorem~\ref{th3.2} above.

\bigskip
A good reference for the material described in this appendix, besides the original papers, is the survey article \cite{CMM} and the book of Yoshino \cite{Y}.

\bigskip

\begin{flushleft}

{\bf AMS Subject Classification: 14M06, 14M07, 13C40, 13C14.}\\[2ex]

%
Robin~HARTSHORNE\\
Department of Mathematics\\
University of California\\
Berkeley, CA\ \ 94720-3840, USA\\
e-mail: \texttt{robin@math.berkeley.edu}\\[2ex]

\end{flushleft}

\end{document}